\documentclass[a4paper,11pt]{article}

\usepackage{amsthm}
\usepackage{caption}
\usepackage{subfigure}
\usepackage{amssymb}
\usepackage{graphicx,subfig}
\usepackage{ebezier}
\usepackage{curves}
\usepackage{epic,eepic}

\newtheorem{theorem}{Theorem}
\newtheorem{lemma}{Lemma}
\newtheorem{proposition}{Proposition}
\newtheorem{definition}{Definition}
\newtheorem{corollary}{Corollary}
\newtheorem{remark}{Remark}

\begin{document}

\title{On the Leaders' Graphical Characterization for Controllability of Path Related Graphs}
\author{Li Dai$^{1*}$, Dianlong Yu $^2$, Zheng Xie $^1$}

\maketitle

\scriptsize{1.College of Liberal Arts and Sciences,
National University of Defense Technology, Changsha 410072, China
2.College  of Intelligent Science and Technology,
 National University of Defense Technology, Changsha 410073, China
 $*$daili@nudt.edu.cn}

\normalsize
\textbf{abstract}
The problem of leaders location plays  an important role in the controllability of undirected graphs.The concept of minimal perfect critical vertex set is introduced by drawing support from the eigenvector of Laplace matrix.
Using the notion of minimal perfect critical vertex set, the problem of finding the minimum number of controllable leader vertices is transformed into the problem of finding all minimal perfect critical vertex sets.
Some necessary and sufficient conditions for special minimal perfect critical vertex sets are provided, such as minimal  perfect critical 2 vertex set, and minimal perfect critical vertex set of path or path related graphs. And further, the leaders location problem for path graphs is solved completely by the algorithm provided in this paper. An interesting result that there never exist a minimal perfect critical 3 vertex set is proved, too.

\textbf{keyword}
controllability, leaders location, multi-agent system, path, generalized star

\section{Introduction}
  Inspired by the swarming behaviors of biological systems and great promises in numerous applications, the field of controllability of multi-agent systems has been studied extensively in recent years \cite{Ali,ZhijianAndHai,ShimaAndMohammad}.By introducing the concept of matching, Liu et al.'s paper \cite{Yangyu} printed in  Nature in 2011 gives a method to find the minimum leaders set for directed networks. However, as pointed out by Ji in \cite{ZhijianAndHai},
when the topological structure of the systems is undirected, how to locate the leaders and what is the minimum number of leaders to insure the controllability are still difficult and largely unknown problems.

 \subsection{Literature Review}
The neighbor-based controllability of undirected graph under a single leader was first formulated by Tanner in  \cite{Tanner} and a necessary and sufficient condition expressed in terms of eigenvalue and eigenvector was derived. In case of multiple leaders, some other algebraic conditions were developed in \cite{Meng,Zhijian&Zidong,Wang&Jiang,Zhijian} etc.
These algebraic conditions lay the foundation for understanding interaction between topological structures of undirected graph and its controllability.
And they are also serve as the theoretical basis of this paper.
The research efforts on characterizing the controllability from a graphical point of view was also motivated by \cite{Tanner} to build controllable topologies. Many kinds of uncontrollable topologies were characterized, such as a symmetric graph with respect to the anchored nodes \cite{RahmaniAndMesbahi}, quotient graphs \cite{MartiniAndEgerstedt} , nodes with the same number of neighbors \cite{JiAndLin}, controllability destructive nodes \cite{Zhijian} etc.
Useful tools and methods  were developed to study the controllability of undirected graph, such as downer branch for tree graphs \cite{ZhijianAndHai}, Zero forcing set\cite{ShimaAndMohammad, Monshizadeh},
 equitable partitions \cite{Meng,Rahmani,Cesar,Martini,Camlibel}, leader and follower subgraphs  \cite{JiAndLin},
$\lambda$-core vertex  \cite{Sciriha,Farrugia}, Distance-to-Leaders (DL) Vector \cite{Yazicioglu}, etc.
Omnicontrollable systems are defined by \cite{Farrugia}, in such systems, the choice of leader vertices that control the follower graph is arbitrary.
Minimal controllability problem (MCP)
that aims to determine the minimum number of state variables that need to be actuated to ensure
system¡¯s controllability was studied in \cite{Olshevsky,Pequito}.
In  study \cite{ZhaoAndGuan}, two
algorithms are established for selecting the fewest leaders to preserve the controllability and the algorithm for leaders¡¯ locations to maximize non-fragility is also designed. Necessary and sufficient conditions to characterize all and only the nodes from which the path or cycle network systerm is controllability were provided in \cite{Liu&Zhijian, parlangeli}.

Although many scholars have devoted themselves to the research in the controllability of undirected graph and achieved many remarkably strong and elegant results, this problem has not been solved yet.
As it is well known that any undirected simple connected graph on $n$ vertices is always $(n-1)$-omnicontrollable. To insure the minimal controllability, which vertices should be selected as leaders is important.
Therefore, Our aim is to find a
method for giving a direct interpretation of the leader vertices  from a graph-theoretic vantage point. In this sense, we provide a new  concept, minimal perfect critical vertex set, to identified the potential leader vertices. This provides a new direction for the study of controllability of undirected systems.

\subsection{Notations and Preliminary Results}
Let $G=(V,E)$ be an undirected and unweighted simple graph, where $V=\{v_1,v_2,\cdots,v_n\}$ is a vertex set and $E=\{v_iv_j|v_i\,\, and \,\,v_j\in V\}$ is an edge set, with an edge $v_iv_j$ is an unordered pair of distinct vertices in $V$. If $v_iv_j\in E$, then $v_i$ and $v_j$ are said to be \textit{adjacent}, or \textit{neighbors}.
$N_S(v_i)=\{v_j\in S| v_iv_j\in E(G)\}$ represents the neighboring set in $S$ of $v_i$, where $S\subset V$. The cardinality of $S$ is denoted by $|S|$. $G[S]$ is the induced subgraph, whose vertex set is $S$ and edge set is $\{v_iv_j\in E(G)|v_i,v_j\in S\}$.
The \textit{valency matrix} $\Delta(G)$ of graph $G$ is a diagonal matrix with rows and columns indexed by $V$, in which the $(i,i)$-entry is the degree of vertex $v_i$, e.g. $|N_G(v_i)|$. Any undirected simple graph  can be represented by its \textit{adjacency matrix}, $D(G)$, which is a symmetric matrix with 0-1 elements. The element in position $(i,j)$ in $D(G)$ is 1 if vertices $v_i$ and $v_j$ are adjacent and 0 otherwise. The symmetric matrix defined as:
\[\mathbf{L}(G)=\mathbf{\Delta}(G)-\mathbf{D}(G)\]
is the \textit{Laplacian } of $G$.The Laplacian is always symmetric and positive semidefinite, and the algebraic multiplicity of its zero eigenvalue is equal to the number of connected components in the graph. For a connected graph, the $n-$dimensional eigenvector associated with the single zero eigenvalue is the vector of ones, $\textbf{1}_n$.

Throughout this paper, it is assumed without loss of generality that $F$ denotes follower vertex set and its vertices  play followers role, and the vertices in $\overline{F}$  are leaders(driver nodes), where $\overline{F}=V \backslash F$ denotes the complement set of $F$.
Let $\mathbf{y}$ be a vector, $\mathbf{y}|_S$ denote the vector obtained from $\mathbf{y}$ after deleting the elements in $\overline{S}$.
Let $\mathbf{L}_{S\rightarrow T}$ denote the matrix obtained from $\mathbf{L}$ after deleting the rows in $\overline{S}$ and columns in $\overline{T}$.
The system described by undirected graph $G$ is said to be controllable (for convenience, $G$ is controllable )if it can be driven from any initial state to any desired state in finite time.
If the followers'  dynamics is (see (4) in \cite{Tanner})
\[\dot{\mathbf{x}}=\mathbf{Ax}+\mathbf{Bu},\]
where $\mathbf{x}$ captures the state of a system which is the stack vector of all  $x_i$ corresponding to follower vertex $v_i\in F$ and $\mathbf{u}$ is the external control
inputs vector which is imposed by the controller and is injected to only some of the vertices, namely the leaders,
the system is controllable with the follower vertex set $F$ if and only if  the
$N \times NM$ controllability matrix
\[\mathbf{C}=[\mathbf{B},\mathbf{AB},\mathbf{A^2B},\cdots,
\mathbf{A^{N-1}B}]\]
has full row rank, that is
$rank(\mathbf{C})=N,$
where $\mathbf{B}=\mathbf{L}_{F\rightarrow \overline{F}}$ and $\mathbf{A}=\mathbf{L}_{F\rightarrow F}$.
This  represents the mathematical condition for controllability, and is well known as Kalman's controllability rank condition\cite{Yangyu, Kalman,Brockett}.

\captionsetup{font={scriptsize}}
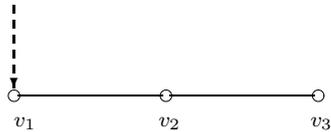
\begin{figure}[h]
\centering
\scriptsize
\setlength{\unitlength}{0.1mm}
\begin{picture}(1000,400)
\put(250,100){\circle{15}}
\put(450,100){\circle{15}}
\put(650,100){\circle{15}}

\put(250,60){$v_1$}
\put(440,60){$v_2$}
\put(640,60){$v_3$}

\put(255,100){\line(1,0){190}}
\put(455,100){\line(1,0){190}}

\put(250,220){\line(0,-1){10}}
\put(250,200){\line(0,-1){10}}
\put(250,180){\line(0,-1){10}}
\put(250,160){\line(0,-1){10}}
\put(250,140){\line(0,-1){10}}
\put(250,120){\vector(0,-1){10}}
\end{picture}
\caption{ The influence of leaders selection on controllability of the systems}
\label{TheLeadersSelectionOn}
\end{figure}

For example, if the vertex $v_1$ is selected as leader, the system is controllable(see fig.\ref{TheLeadersSelectionOn}).But, if $v_2$ plays the leaders role, it is NOT controllable. This paper will address the graphical characterization of leaders to insure the systems's controllability.

In most real systems such as multi-agent systems or complex networks, we are particularly interested in identifying the minimum number of leaders, whose control is sufficient and fully control the system's dynamics.

 In term of eigenvalues and eigenvectors of submatrices of Laplace, \cite{Meng,Zhijian&Zidong,Wang&Jiang,Zhijian} presented a necessary and sufficient algebraic condition on controllability.

\begin{proposition}\label{proposition1}
\cite{Zhijian&Zidong,Wang&Jiang,Zhijian}
The undirected graph $G$ is controllable under the leader vertex set $\overline{F}$ if and only if $\mathbf{L}$ and $\mathbf{L}_{F\rightarrow F}$  share no common eigenvalue.
\end{proposition}

\begin{proposition}\label{proposition2}
\cite{Meng,Zhijian}
The undirected graph $G$ is controllable under the leader vertex set $\overline{F} $ if and only if $\mathbf{y}|_{\overline{F}}\neq \mathbf{0}$ \, (\,\,$\forall\,\,\mathbf{y}$  a eigenvector of $\mathbf{L}$).
\end{proposition}

 The Proposition  \ref{proposition2} gives the algebraic characteristics of leader vertex set. It is worth noting that the eigenvector $\textbf{y}$ in Proposition \ref{proposition2} has the characteristic of arbitrary. Therefore, when $\textbf{L}$ has multiple eigenvalues, it is not possible to draw a conclusion only by examining all the linearly independent eigenvectors, but also by further verifying all the eigenvectors with zero components. From the point of view of numerical calculation, this verification is too computational and difficult to implement. It is clear that the topology of the interconnection graph  $G$ completely determines its controllability properties. So, this paper will focus on the graph theoretic characterization of the leader vertices.

The remainder of this paper is organized as follows. In
Section \ref{section2}, we provide three new  concepts: critical vertex set, perfect critical vertex set and minimal perfect critical vertex set. Necessary and sufficient conditions  for $S$ to be a minimal perfect critical 2 vertex set is presented. An interesting result that there never exist minimal perfect critical 3 vertex set is also proved in Section \ref{section2}.
Section \ref{section3} is the main part of this paper. In this section, we provide a algorithm to locate all leader vertices of path by finding out its all minimal perfect critical vertex set.
Graphs constructed by adding paths incident to one vertex $v_0$ are investiaged in Section \ref{section4}.
Finally, our conclusions are summarized  in Section \ref{section5}.

\section{Minimal Perfect Critical Vertex Set}\label{section2}
According to proposition \ref{proposition2}, for any $S\subset V$ and $S\neq \emptyset$,
if there exist an eigenvector $\textbf{y}$ of Laplace matrix $\textbf{L}$ such that $\mathbf{y}|_{S}=\mathbf{0}$, then $S$ cannot be used as a leader vertex set.
So, in order to locate the leaders   of graph $G$, the following concepts are proposed.

\subsection{Three Definitions}\label{subsection2.1}

\begin{definition}\label{definition1}
(\textbf{critical vertex set})\,\,Let $S$ be a nonempty subset of $V$, if there exist an eigenvector $\textbf{y}$ such that $\textbf{y}|_{\overline{S}}=\textbf{0}$, then $S$ is called a critical vertex set(CVS) and $\textbf{y}$ is a inducing eigenvector.  $S$ is called a critical $k$ vertex set, if $|S|=k$.
\end{definition}

\begin{definition}\label{definition2}
(\textbf{perfect critical vertex set})\,\,
Let $S$ be a critical vertex set, if there exist a eigenvector $\textbf{y}$ satisfy that $\textbf{y}|_{\{v_i\}}\ne 0 (\forall v_i\in S)$, then $S$ is called a perfect critical vertex set(PCVS). And $S$ is called a perfect critical $k$ vertex set, if $|S|=k$.
\end{definition}

\begin{definition}\label{definition 3}
(\textbf{minimal perfect critical vertex set})
A perfect critical vertex set is called a minimal perfect critical vertex set (MPCVS) if its any proper subset is no longer a perfect critical vertex set. And $S$ is called a minimal perfect critical $k$ vertex set, if $|S|=k$.
\end{definition}

\begin{remark}\label{remark-definition}
By the Definitions, $V$ is a trivially CVS and a PCVS induced by the eigenvector $\textbf{1}_n$. $V$ is a MPCVS if and only if $G$ is controllable under any single vertex selected as leader, e.g. $G$ is omnicontrollable.
\end{remark}

\begin{figure}[h]
\centering
\scriptsize
\setlength{\unitlength}{1mm}
\begin{picture}(40,40)
\put(0,5){\circle{1.5}}
\put(30,5){\circle{1.5}}
\put(0,20){\circle{1.5}}
\put(0,34.5){\circle{1.5}}
\put(15,20){\circle{1.5}}
\put(43,18.5){\circle{1.5}}
\put(27.8,33.2){\circle{1.5}}

\put(-5,4){$v_1$}
\put(-5,19){$v_2$}
\put(-5,34){$v_3$}
\put(13,23){$v_4$}
\put(32,4){$v_7$}
\put(42.5,20){$v_6$}
\put(28,35){$v_5$}

\put(0.6,5.5){\line(1,1){14}}
\put(15.4,20.6){\line(1,1){12}}%v4-v5
\put(0.6,20){\line(1,0){13.8}}%v2-v4
\put(14.5,20){\line(-1,1){14}}%v4--v3
\put(42.5,19){\line(-1,1){14}}%v6--v5
\put(29.5,5.5){\line(-1,1){14}}%v7--v4
\put(30.5,5.5){\line(1,1){12.3}}%v7--v4

\end{picture}

\caption{Graph $G$ with 4 distinct minimal perfect critical vertex set}
\label{figureForDifferentMPCVS}%ͼ2
\end{figure}
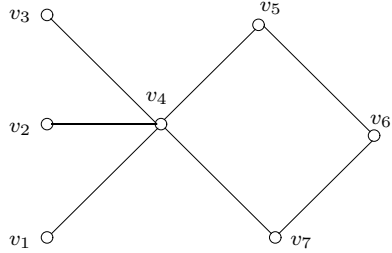

For example, see fig.\ref{figureForDifferentMPCVS}, $\{v_1,v_5\}$ is not a CVS, $\{v_1,v_2,v_3,v_4\}$ is a CVS but not a PCVS. $S=\{v_1,v_2,v_3\}$ is a PCVS but not a MPCVS. $S_1=\{v_1,v_2\}$, $S_2=\{v_1,v_3\}$, $S_3=\{v_2,v_3\}$ and $S_4=\{v_5,v_7\}$ are all MPCVSs of the graph in  fig.\ref{figureForDifferentMPCVS} .

\begin{remark}\label{remarkForPropositionStatedWithMPVCS}
 Proposition \ref{proposition2} can be restated as: The undirected graph $G$ is controllable under the leader vertex set $\overline{F} $ if and only if for each MPCVS $S$ , $S \bigcap \overline{F}\ne \emptyset$.
\end{remark}

It is Remark \ref{remarkForPropositionStatedWithMPVCS} that inspired us to study MPCVS. Because, from  Remark \ref{remarkForPropositionStatedWithMPVCS} , there is a close relationship between MPCVS and the minimum leaders set.
In other words, when we find all MPCVS of $G$, we find the minimum leader set and hence the minimum number of leader vertices. For example, by Remark \ref{remarkForPropositionStatedWithMPVCS} and all its 4 MPCVS above, graph $G$ in fig.\ref{figureForDifferentMPCVS} is not controllable under any single leader because $\bigcap_{i=1}^{4}S_i=\emptyset$ , or any two vertices. So, the minimum leaders set are $\{v_i,v_j,v_k\}$, where $v_i,v_j$ comes from $\{v_1,v_2,v_3\}$ and $v_k$ from $\{v_5,v_7\}$. Therefor, the minmum number of leaders is 3.

Moreover, many MPCVSs have typical graphical characteristics. For example, all of the 4 MPCVSs of $G$ in fig.\ref{figureForDifferentMPCVS} have the graphical structure stated in  Theorem \ref{theorem3}. This is another reason for us to investigate MPCVS.

\subsection{Sufficient Conditions for Critical Vertex Set}\label{subsection2.2}
For undirected graph, Laplacian matrix $\textbf{L}$ is  symmetric, all the eigenvectors are orthogonal to each other, so knowing $\textbf{1}_n$ is an eigenvector of $\textbf{L}$, it is immediate that all the other eigenvectors of $\textbf{L}$ are orthogonal to $\bold{1}_n$, that is, for all eigenvector $ \bold{y}$,

\begin{equation}\label{eq1}
\textbf{1}_{n}^T\textbf{y}=\sum_{i=1}^{n}y_i=0.
\end{equation}

The equality in (\ref{eq1}) is useful throughout the paper.

If $S$ is a CVS, then
\begin{equation}\label{equationFor|S|>=2}
|S|\geq 2.
\end{equation}
 In fact,
Suppose $|S|=1$, without loss of generality, $S=\{v_1\}$.
Let $\textbf{y}=(y_1,y_2,\cdots,y_n)^T$ be the inducing eigenvector  associated with eigenvalue $\lambda$, then $\textbf{Ly}=\lambda \textbf{y}$ and $\textbf{y}|_{\overline{S}}=\textbf{0}$.
By $\textbf{y}|_{\overline{S}}=\textbf{0}$ and (\ref{eq1}), $\textbf{y}|_{S}=\textbf{0}$, e.g. $\textbf{y}=\textbf{0}$, this is
in contradiction with the fact that $\textbf{y}$ is an eigenvector.
Further, since any subset $S$ with $|S|=1$ isn't a CVS, by Remark \ref{remarkForPropositionStatedWithMPVCS}, $G$  is controllable with the leader vertex set $\overline{F}$ when $|F|=1$.

Now, we are going to investigate the properties of critical vertex set. Firstly, a sufficient conditions for $S$ to be a CVS is provided in the following Proposition \ref{proposition3}, which describes a special case of the symmetry-based uncontrollability results.

\begin{proposition}\label{proposition3}
Let $G$ be an undirected connected graph of order $n$, $S\subset V$ and $|S|\geq 2$, if for any $v\in \overline{S}$, either $N_S(v)=\emptyset$ or $N_S(v)=S$, then $S$ is a critical vertex set.
\end{proposition}

\textbf{Proof}\quad Let $|\{v\in \overline{S}|N_S(v)=S\}|=m$, then
\[\textbf{L}_{S\rightarrow S}-m\textbf{I}_{|S|}\]
is Laplacian of subgraph $G[S]$, where $\textbf{I}_{|S|}$ denotes the $|S|$ dimensional identity matrix. Considered (\ref{eq1}), there exist an eigenvector $\textbf{y}_S$ of the Laplacian $\textbf{L}_{S\rightarrow S}-m\textbf{I}_{|S|}$ such that $\textbf{1}_{|S|}^T\textbf{y}_S=0$.

Set vector $\textbf{y}$ as $\textbf{y}|_S=\textbf{y}_S$ and $\textbf{y}|_{\overline{S}}=\textbf{0}$. It can be seen that

 \begin{equation}\label{eq2}
 \textbf{Ly}=\left[ \begin{array}{ll}

  \textbf{L}_{S\rightarrow S}& \textbf{L}_{S\rightarrow \overline{S}} \\
  \textbf{L}_{\overline{S}\rightarrow S}& \textbf{L}_{\overline{S}\rightarrow \overline{S}}
  \end{array}\right]
  \left[\begin{array}{c} \textbf{y}_S\\ \textbf{0 }\end{array}\right]=\left[\begin{array}{c}  \lambda \textbf{y}_S\\  \textbf{L}_{\overline{S}\rightarrow S}\textbf{y}_S \end{array}\right].
 \end{equation}
Noticing that the rows in matrix $\textbf{L}_{\overline{S}\rightarrow S}$ are either  ones or zeros, the conclusion is proved by $\textbf{1}_{|S|}^T\textbf{y}_S=0$.   \qed

For example, by Proposition \ref{proposition3}, $\{v_1,v_3\}$ is a CVS, therefor, the graph $G$ in fig.\ref{TheLeadersSelectionOn} is uncontrollable when $v_2$ is selected as leader.

\begin{remark}\label{remark3}
The condition provided in Proposition \ref{proposition3}  implies that some critical vertex sets are closely related to equitable partitions. For example, let $S$ be the followers and $\overline{S}$ be the leaders. From earlier results in the literature \cite{Martini}, we know that in the case of Proposition \ref{proposition3} the maximal relaxed equitable partition would put all the leaders into a single cell, hence the system is uncontrollable. But, some other critical vertex sets have nothing to do with equitable partitions or almost equitable partitions(AEP, see \cite{Cesar}). For example, let $S$ be the perfect critical vertex set in fig.\ref{PerfectCriticalVertexSetIsCloselyRelatedTo}(a). The partitions obtained by putting all the vertices in $\overline{S}$ into a single cell is not a AEP.
\end{remark}

\subsection{Minimal Perfect Critical 2  and  3 Vertex Set }
Armed with the above properties, critical $k$ vertex set with $k\leq 3$ can be determined directly from their graphical characterization.
This is achieved via a detailed analysis of the inducing eigenvector.
\begin{lemma}\label{lemma1}
Let $G$ be an undirected connected graph and $S$ be a perfect critical $k$ vertices set, then for any $v\in \overline{S}$, $|N_S(v)|\ne 1$ and $|N_S(v)|\ne k-1$.
\end{lemma}
\textbf{Proof}\quad Let $S=\{v_1,v_2,\cdots,v_k\}$ be a perfect critical vertices set and $\textbf{y}=(y_1,y_2,\cdots,y_k,0,0,\cdots,0)^T$ be the inducing eigenvector. $y_i\ne 0(\forall 1\leq i \leq k)$ since $S$ is a perfect critical vertices set.

$\forall v\in \overline{S}$, suppose the $|N_S(v)|=1$, without loss of generality, say, $vv_1\in E$ and $vv_i\notin E(G)(\forall i\neq 1)$, then $\textbf{L}_{\{v\}\rightarrow V}\textbf{y}=y_1\neq 0$. On the other hand, $\textbf{y}|_{\overline{S}}=\textbf{0}$ and $v \in \overline{S}$, so
$\textbf{L}_{\{v\}\rightarrow V}\textbf{y}=0$, this is a contradiction.

Together with (\ref{eq1}), $|N_S(v)|\ne k-1$ can be proved similarly.\qed

By (\ref{equationFor|S|>=2}), critical 2 vertex set is also a minimal perfect critical 2 vertex set. The following Theorem \ref{theorem3} will follow from Lemma \ref{lemma1} and Proposition \ref{proposition3}.

\begin{theorem}\label{theorem3}
Let $G$ be an undirected connected graph, $S\subset V$ and $|S|=2$, then $S$ is a minimal  perfect critical 2 vertex set if and only if $\forall v\in \overline{S}$, either $N_S(v)=\emptyset$ or $N_S(v)=S$.\qed
\end{theorem}

For example, see graph $G$ in fig.\ref{figureForDifferentMPCVS},
all its 4 MPCVSs  can be recognized by the graphical characterization stated in Theorem \ref{theorem3}.

\begin{remark}\label{remark5}
 FromTheorem \ref{theorem3}, we know that the perfect critical 2 vertex set is what named \textit{twins nodes} by \cite{Biyikoglu} and also it is double controllability destructive (DCD) node tuple given by \cite{Zhijian} and \cite{Liu&Zhijian}. Hence, one can see that the perfect critical vertex set is the extension and generalization of twins nodes and controllability destructive nodes.
\end{remark}

 But the minimal perfect critical 3 vertex set is not the same as the triple controllability destructive nodes  \textit{TCD nodes} named by \cite{Zhijian}, because we will prove that there do not exist a minimal perfect critical 3 vertex set. That is the following Theorem \ref{theorem4}.

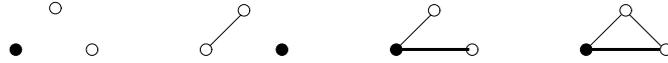
\begin{figure}[h]
\centering
\setlength{\unitlength}{1mm}
\begin{picture}(80,20)
\put(0,5){\circle*{1.5}}
\put(10,5){\circle{1.5}}
\put(25,5){\circle{1.5}}
\put(35,5){\circle*{1.5}}
\put(50,5){\circle*{1.5}}
\put(60,5){\circle{1.5}}
\put(75,5){\circle*{1.5}}
\put(85.5,5){\circle{1.5}}

\put(25.5,5.5){\line(1,1){4.1}}
\put(5.2,10.5){\circle{1.5}}
\put(30,10.2){\circle{1.5}}
\put(50.5,5.5){\line(1,1){4.1}}
\put(55,10.2){\circle{1.5}}
\put(50.5,5){\line(1,0){9}}
\put(75.5,5.5){\line(1,1){4.1}}
\put(80.2,10.2){\circle{1.5}}
\put(80.7,9.6){\line(1,-1){4.1}}
\put(75.7,5){\line(1,0){9}}

\end{picture}
\caption{all possible topology structures of $G[S]$ with $|S|=3$}
\label{AllPossibleTopologyStructures}
\end{figure}

\begin{theorem}\label{theorem4}
Let $G$ be an undirected connected graph, $S\subset V$ and $|S|=3$, then $S$ is NOT a minimal perfect critical vertex set.
\end{theorem}
\textbf{Proof}\quad Suppose $S$ is a minimal perfect critical vertex set.Consider the subgraph $G[S]$, all  4  possible topology structures of $G[S]$ are depicted in fig. \ref{AllPossibleTopologyStructures}.

%%%%%%%%%%%  ͼ 3      %%%%%%%%%%%
%%%%%%%%%%%% ͼ 3  ½á Êø    %%%%%%%%%

For each topology of $G[S]$ in fig.\ref{AllPossibleTopologyStructures}, let $T$ be the vertex set of white vertices and $u$ be the black vertex, one can have either $N_T(u)=\emptyset$ or $N_T(u)=T$.

 By Lemma \ref{lemma1}, $\forall v\in\overline{S}$, either $|N_S(v)|=0$ or $|N_S(v)|=3$.

Noticing that $\overline{T}=\{u\}\bigcup\overline{S}$, by Proposition \ref{proposition3}, $T$ is a critical vertex set. This contradicts the assumption.\qed

\begin{remark}\label{remark4}
There do exist some perfect critical 3 vertex set, but there do not exist any minimal perfect critical 3 vertex set. For example, see fig.\ref{figureForDifferentMPCVS} , $S=\{v_1,v_2,v_3\}$ is a perfect critical vertex set because there exist a eigenvector $\textbf{y}$ such that $\textbf{y}|_{\overline{S}}=\textbf{0}$ and $\textbf{y}|_{v_i}\neq 0(\forall v_i\in S)$. But, by Theorem \ref{theorem3} and Definition \ref{definition 3}, $S$ is not a MPCVS.
\end{remark}

%%%%%%%%%%%  ͼ 4      %%%%%%%%%%%
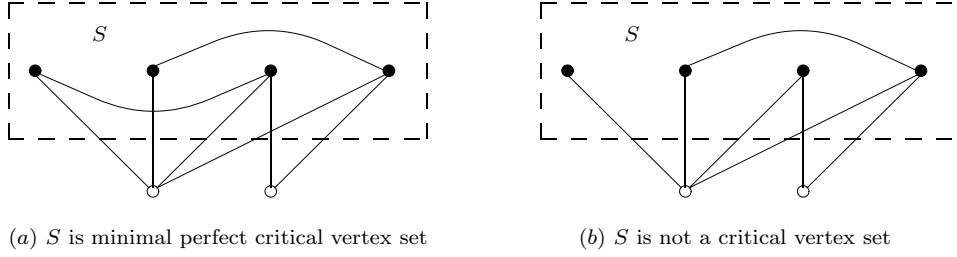
\begin{figure}[h]
\centering
\setlength{\unitlength}{1mm}
\begin{picture}(120,40)
\put(20,5){\circle{1.5}}
\put(20,5.5){\line(0,1){15}}
\put(20,21){\circle*{1.5}}
\put(19.5,5.5){\line(-1,1){15}}
\put(4.5,21){\circle*{1.5}}
\put(20.5,5.5){\line(1,1){15}}
\put(35.5,21){\circle*{1.5}}
\put(35.5,5){\circle{1.5}}
\put(35.5,5.5){\line(0,1){15}}
\put(36,5.5){\line(1,1){15}}
\put(51,21){\circle*{1.5}}
\put(20.5,5.5){\line(2,1){30}}
\put(1,12){\dashbox{2}(55,18)}
\spline(5.2,20.5)(20,14)(35,20.5)
\spline(20,21.5)(35.5,28)(50.5,21)
\put(12,25){\scriptsize $S$}
\put(1,-2){\scriptsize  $(a)\,\, S$ is minimal perfect critical vertex set}

\put(90,5){\circle{1.5}}
\put(90,5.5){\line(0,1){15}}
\put(90,21){\circle*{1.5}}
\put(89.5,5.5){\line(-1,1){15}}
\put(74.5,21){\circle*{1.5}}
\put(90.5,5.5){\line(1,1){15}}
\put(105.5,21){\circle*{1.5}}
\put(105.5,5){\circle{1.5}}
\put(105.5,5.5){\line(0,1){15}}
\put(106,5.5){\line(1,1){15}}
\put(121,21){\circle*{1.5}}
\put(90.5,5.5){\line(2,1){30}}

\put(71,12){\dashbox{2}(55,18)}
\spline(90,21.5)(105.5,28)(120.5,21)
\put(82,25){\scriptsize$S$}
\put(76,-2){\scriptsize $(b)\,\, S$ is not a critical vertex set}
\end{picture}
\vspace{0.5cm}
\caption{all possible topology structures of $G[S]$ with $|S|=3$}
\label{PerfectCriticalVertexSetIsCloselyRelatedTo}
\end{figure}
%%%%%%%%%%%%%%  ͼ4 ½áÊø %%%%%%%%%%%%

Although there does not exist a minimal perfect critical 3 vertex set,  minimal perfect critical 4 vertex set does exist, see fig.\ref{PerfectCriticalVertexSetIsCloselyRelatedTo}(a).
From Theorem \ref{theorem3}, perfect critical 2 vertex set is completely determined by the relationship between $\overline{S}$ and $S$, and have nothing to do with the interconnection topology of subgraph $G[S]$. But, unlike the perfect critical 2 vertex set,
the topology structure of $G[S]$ will have an effect on whether $S$ is a minimal perfect critical 4 vertex set or not, see fig.\ref{PerfectCriticalVertexSetIsCloselyRelatedTo}.
The virtue that perfect critical $k(k\geq 4)$ vertex set should have was needed to be characterized from both algebraic and graphical perspectives.
Developing such a characterization is along the directions of our current research.

\section{Minimal Perfect Critical Vertex Set of Path}\label{section3}
In this section, we will solve the leaders location problem for path completely by means of MPCVS.

\subsection{Spectral Propoerties}\label{subsection3.1}

A path graph $P_n$ is a finite sequence of vertex $v_1,v_2,\cdots,v_n$ starting with $v_1$ and ending with $v_n$ such that consecutive vertex are adjacent. A subset $S\subset V$  is said to be isolated vertex set where there are no edges among the verties in $S$.

 If $S$ be a perfect critical vertex set of path $P_n$, by Lemma \ref{lemma1},  $\overline{S}$ must be isolated vertex set. So, without loss of generality, let $\overline{S}=\{v_{i_1},v_{i_2},\cdots,v_{i_k}\}$ be a isolated vertex set and $1< i_1< i_2<\cdots < i_k < n$.
Let $S_{i_0}=\{v_1,v_2,\cdots,v_{i_1-1}\}$,
$S_{i_1}=\{v_{i_1+1},v_{i_1+2},\cdots,v_{i_2-1}\},\cdots$,
$S_{i_k}=\{v_{i_k+1},v_{i_k+2},\cdots,v_{i_n}\}$.
Recall Lemma \ref{lemma1}, we know that $1<i_1$ and $i_k<n$, e.g.,

\begin{equation}\label{m>1}
|S_{i_0}|\geq 1\,\, and \,\,|S_{i_k}|\geq 1.
\end{equation}

It is easy to see that the matrix $\textbf{L}_{S\rightarrow S}$ is a block matrix with the following form

\begin{equation}\label{LStoS=blocked}
\textbf{L}_{S\rightarrow S}=\left[\begin{array}{cccc}
 \textbf{L}_{S_{i_0}\rightarrow S_{i_0}}  & \textbf{0}  & \textbf{0}  &\textbf{0}  \\
   \textbf{0} &  \textbf{L}_{S_{i_1}\rightarrow S_{i_1}}  & \textbf{0}  &  \textbf{0} \\
 \textbf{0}   & \textbf{0} & \ddots  & \textbf{0}  \\
   \textbf{0} & \textbf{0}  & \textbf{0}  &  \textbf{L}_{S_{i_k}\rightarrow S_{i_k}}
\end{array}\right]
\end{equation}

By rearrange the columns and rows we can always write the Laplacian of $P_n$ into the following form:

\begin{equation}\label{L=blocked}
\textbf{L}=\left[
\begin{array}{cc}
\textbf{L}_{S\rightarrow S}  &   \textbf{L}_{S\rightarrow \overline{S}}\\
\textbf{L}_{\overline{S}\rightarrow S}  &  \textbf{L}_{\overline{S}\rightarrow\overline{S}}
\end{array}
\right].
\end{equation}

Since for path $P_n$, only the consecutive vertices are adjacent, so there are exactly 2 elements are -1 while the other elements are 0 for every  vector
$\textbf{L}_{\{v_{i_j}\}\rightarrow S}(j=1,2,\cdots,k)$. That is

\begin{equation}\label{Lvi-S}
\textbf{L}_{\{v_{i_j}\}\rightarrow S}=(0,0,\cdots,0,
\stackrel{\stackrel{v_{i_{j-1}}}{\uparrow}}{-1},
\stackrel{\stackrel{v_{i_{j+1}}}{\uparrow}}{-1},0,\cdots,0)
(j=1,2,\cdots,k)
\end{equation}

For path, the matrixes in right side of (\ref{LStoS=blocked}) have similar structure. Let $m$ and $M$ be the dimension of  the following useful matrix $D_m$ and $B_M$, respectively. These matrices  play an important role to determine the locations of leaders under which the controllability of paths  can be realized. Thus,
the first submatrix in (\ref{LStoS=blocked}) can be written as $\textbf{L}_{S_{i_0}\rightarrow S_{i_0}}=\textbf{D}_{|S_{i_0}|}$. By symmetric permutation reversing all the components, the last submatrix in (\ref{LStoS=blocked}) can be written as $\textbf{L}_{S_{i_k}\rightarrow S_{i_k}}=\textbf{D}_{|S_{i_k}|}$. The other submatrices in (\ref{LStoS=blocked}) can be written as $\textbf{L}_{S_{i_j}\rightarrow S_{i_j}}=\textbf{B}_{|S_{i_j}|}(j=1,2,\cdots,k-1)$.

\[\textbf{D}_m=\left[
\begin{array}{cccccccc}
1 & -1 & 0  &  0  & \cdots & 0 & 0 & 0\\
-1 &  2  & -1 & 0  &  \cdots  & 0 &  0  &0 \\
\vdots &  \vdots &  \vdots&  \vdots &    & \vdots  & \vdots &  \vdots  \\
0&  0 &  0 &  0 &  \cdots  & 0  & -1 &  2 \\
\end{array}
\right]_{m\times m},\]

\[\textbf{B}_M=\left[
\begin{array}{cccccccc}
2 & -1 & 0  &  0  & \cdots & 0 & 0 & 0\\
-1 &  2  & -1 & 0  &  \cdots  & 0 &  0  &0 \\
\vdots &  \vdots &  \vdots&  \vdots &    & \vdots  & \vdots &  \vdots  \\
0&  0 &  0 &  0 &  \cdots  & 0  & -1 &  2 \\
\end{array}
\right]_{M\times M}.\]

Naturally, We are going to investigate the spectral properties of $\textbf{D}_m$ and $\textbf{B}_M$.

For convenient,  we introduce some useful notations. For any $\lambda\in[0,4]$, $\theta$ is called a angle associated with $\lambda$, where $\theta$ is defined as $\cos\theta=\frac{2-\lambda}{2}$ , $\sin\theta=\frac{\sqrt{4\lambda-\lambda^2}}{2}$ and $\theta\in[0,\pi]$. If $\lambda$ is a eigenvalue, then $\theta$ is called a eigenangle.
Let $\phi_i(\lambda)$  and $\psi_i(\lambda)$ be the $i$-th sequential principal minor of $\det(\lambda \textbf{I}-\textbf{D}_m)$ and $\det(\lambda \textbf{I}-\textbf{B}_M)$, respectively, then  we have the following useful lemmas.

\begin{proposition}\label{spectralOf-Dm}
Let $\textbf{y}=(y_1,y_2,\cdots,y_m)$ be an eigenvector of $\textbf{D}_m$.

(i)\quad $\lambda$ is a eigenvalue of $\textbf{D}_m$ if and only if
$\theta\in\left\{\frac{(2l-1)\pi}{2m+1}|1\leq l \leq m\right\}$ is associated with $\lambda$.

(ii)\quad $y_i=\phi_{i-1}(\lambda)y_1,(i=2,3,\cdots,m).$
\end{proposition}

\textbf{Proof}\quad (i) By applying the Laplace expansion to the last row of $\lambda \textbf{I}-\textbf{D}_m$, the following recurrence formula hold:
\[\phi_m(\lambda)=(2-\lambda)\phi_{m-1}(\lambda)-\phi_{m-2}(\lambda).\]
Then, from Gersgorin disk theorem, it follows that the eigenvalue $\lambda\leq 4$, that means $\lambda^2-4\lambda\leq 0$. Solving this recurrence and taking $\phi_1(\lambda)=1-\lambda$ into consideration, we have

\begin{equation}\label{theta-of-Dm=0}
\phi_m(\lambda)=\frac{\cos{\frac{(2m+1)\theta}{2}}}{\cos\frac{\theta}{2}}.
\end{equation}
Thus concluding the first part of the proof.

(ii) The claim can be verified via mathematical induction from the fact:
$(\lambda\textbf{I}-\textbf{D}_m)\textbf{y}=\textbf{0}.$  \qed

Similarly, we have

\begin{proposition}\label{spectralOf-BM}
Let $\textbf{y}=(y_1,y_2,\cdots,y_M)$ be an eigenvector of $\textbf{B}_M$.

(i)\quad $\lambda$ is a eigenvalue of $\textbf{B}_M$ if and only if
$\theta\in\left\{\frac{h\pi}{M+1}|1\leq h \leq M\right\}$ is associated with $\lambda$.

(ii)\quad $y_i=\psi_{i-1}(\lambda)y_1,(i=2,3,\cdots,M).$
\end{proposition}

\begin{lemma}\label{common-eigenvalue-of-Lsiktosik}
If $S$ is a perfect critical vertex set of $P_n$, then  all $\textbf{L}_{S_{i_j}\rightarrow S_{i_j}}(j=0,1,\cdots,k)$ in (\ref{LStoS=blocked}) have at least one common eigenvalue.
\end{lemma}

\textbf{Proof}\quad If $S$ is a perfect critical vertex set , then there exist an eigenvector $\textbf{y}$ such that $\textbf{y}|_{\overline{S}}=\textbf{0}$ and $\textbf{y}|_{v\in S}\neq 0$.
From (\ref{L=blocked}), we have
$\textbf{L}_{S\rightarrow S}\textbf{y}|_S=\lambda\textbf{y}|_S.$
Now, consider (\ref{LStoS=blocked}), all of $\textbf{L}_{S_{i_j}\rightarrow S_{i_j}}(j=0,1,\cdots,k)$ have at least one common eigenvalue $\lambda$ because $\textbf{y}|_{S_{i_j}}\neq \textbf{0}$ are the corresponding eigenvectors.\qed

\begin{lemma}\label{equal-dimension-for-BM}
If $S$ is a perfect critical vertex set of $P_n$, for all $S_{i_j}(j=0,1,\cdots,k)$ in (\ref{LStoS=blocked}), the following equalities holds:

(i)\quad  $|S_{i_0}|=|S_{i_k}|.$

(ii)\quad  $|S_{i_1}|=|S_{i_2}|=\cdots=|S_{i_{k-1}}|.$

\end{lemma}

\textbf{Proof}\quad (i) Set $|S_{i_0}|=m_1$ and $|S_{i_k}|=m_2$. Suppose that $m_1\neq m_2$. Without loss of generality, $m_2>m_1$. From  Lemma \ref{common-eigenvalue-of-Lsiktosik},  $\textbf{L}_{S_{i_0}\rightarrow S_{i_0}}(=\textbf{D}_{m_1})$ and $\textbf{L}_{S_{i_k}\rightarrow S_{i_k}}(=\textbf{D}_{m_2})$ share a eigenvalue $\tilde{\lambda}$, together with Proposition \ref{spectralOf-Dm}(i), they share eigenangle $\tilde{\theta}$. That means
there exist $l_1,l_2$ such that
\[\frac{2l_1-1}{2m_1+1}\pi=\frac{2l_2-1}{2m_2+1}\pi=\tilde{\theta}.\]
Since $m_2>m_1$ and $\phi_{m_1}(\tilde{\lambda})$ is the $m_1$-th sequential principal minor of $\phi_{m_2}(\tilde{\lambda})$, by Proposition \ref{spectralOf-Dm}(ii), the $(m_1+1)$-th element of the eigenvector $\textbf{y}|_{S_{i_k}}$ is zero. This contradicts the fact that $S$ is a perfect critical vertex set.

The proof of (ii) can be carried out in the same manner as (i).\qed

\begin{lemma}\label{|Dm|=|BM|}
If $S$ is a perfect critical vertex set of $P_n$, for vertex set $S_{i_0}$ and $S_{i_1}$ in (\ref{LStoS=blocked}), then either $S_{i_1}$ is  empty or $|S_{i_1}|=2|S_{i_0}|.$
\end{lemma}

\textbf{Proof}\quad If $S_{i_1}$ is empty, the proof is trivial; thus, let $|S_{i_1}|=M>0$ and $|S_{i_0}|=m$.

From Lemma \ref{common-eigenvalue-of-Lsiktosik}, $\textbf{L}_{S_{i_0}\rightarrow S_{i_0}}$ and $\textbf{L}_{S_{i_1}\rightarrow S_{i_1}}$ have a common eigenvalue $\tilde{\lambda}$. Let $\textbf{y}$ be the induced eigenvector of the  perfect critical vertex set $S$.
According to Proposition \ref{spectralOf-Dm}(i) and Proposition \ref{spectralOf-BM}(i), there exist numbers $l,h$, for $1\leq l_0\leq m, 1\leq h_0\leq M$, such that
\begin{equation}\label{M=2m}
\frac{(2l_0-1)}{2m+1}=\frac{h_0}{M+1}.
\end{equation}
We claim that $2l_0-1$ and $2m+1$ are coprime. Otherwise, recall the proof of Lemma \ref{equal-dimension-for-BM}, we know that there exist at least one entry of the eigenvector $\textbf{y}|_{S_{i_0}}$ vanish.  Similarly, $h_0$ and $M+1$ are coprime, too. That implies $2m+1=M+1$, e.g. $M=2m$.  \qed

\begin{lemma}\label{lemma5}
If $S$ is a minimal perfect critical vertex set, for vertex set $S_{i_0}$ in (\ref{LStoS=blocked}), let $|S_{i_0}|=m$, then $2m+1$ is an odd prime.
\end{lemma}
\textbf{Proof}\quad From Lemma \ref{|Dm|=|BM|},  Proposition \ref{spectralOf-Dm}(i) and Proposition \ref{spectralOf-BM}(i),  we know that the submatrices $\textbf{L}_{S_{i_j}\rightarrow S_{i_j}}$ have the following $m$ common eigenangles:

\begin{equation}\label{eq3}
\{\frac{1}{2m+1}\pi,\frac{3}{2m+1}\pi,\cdots,\frac{2m-1}{2m+1}\pi\}.
\end{equation}

Suppose $2m+1$ is not a prime, there exist two factors $p_1,p_2$  such that $2m+1=p_1p_2$. Since $2m+1$ is odd, both $p_1$ and $p_2$ are odd.
Therefor, $\frac{1}{p_2}\pi$ is one of the eigenangle in (\ref{eq3}).
By Proposition \ref{spectralOf-Dm}(ii), the $(p_2+1)$-th element of the eigenvector associated with eigenvalue $\frac{p_1}{2m+1}\pi$ is zero. This is a contradiction to the fact that $S$ is a MPCVS.  \qed

\subsection{Equivalence Characterization of MPCVSs of Path Graphs}

The following Theorem \ref{theVertexNumberOfPath} provided a equivalence characterization of MPCVS of path graph.

\begin{theorem}\label{theVertexNumberOfPath}
Let $S$ be a vertex set of path $P_n$ and $\overline{S}=\{v_{i_1},v_{i_2},\cdots,v_{i_k}\}$ is isolated. Let $|S_{i_0}|=m$, $|S_{i_1}|=M$.
 Then $S$ is a minimal perfect critical vertex set if and only if the following assertions hold:

 (i)  $|S_{i_k}|=|S_{i_0}|$, $|S_{i_1}|=|S_{i_2}|=\cdots=|S_{i_{k-1}}|$.

 (ii)  $M=0$ or $M=2m$.

 (iii)  $2m+1$ is a odd prime.
\end{theorem}

\textbf{Proof}\quad The necessity is proved in Lemma \ref{equal-dimension-for-BM}, Lemma \ref{|Dm|=|BM|} and Lemma \ref{lemma5}.

Sufficiency: Case 1 $M>0$.

 Since $M=2m$,  all of the eigenangles in (\ref{eq3}) are common eigenangles of all submatrices $\textbf{L}_{S_{i_i}\rightarrow S_{i_j}}(j=0,1,\cdots,k)$.
 None of the eigenangles in (\ref{eq3}) is a eigenangle of any sequential principal minor of $\textbf{L}_{S_{i_j}\rightarrow S_{i_j}}(j=0,1,\cdots,k)$ because of the condition (iii). So,  from  Proposition \ref{spectralOf-Dm}(ii) and Proposition \ref{spectralOf-BM}(ii), we know any eigenvectors $\textbf{L}_{S_{i_j}\rightarrow S_{i_j}}(j=0,1,\cdots,k)$ associated with the common eigenangles  have zero elements.
 Therefor, we only need to proof that there exist a eigenvector $\textbf{y}$ of $\textbf{L}$ such that $\textbf{y}|_{\overline{S}}=\textbf{0}$.

Arbitrarily selecta common eigenangle $\theta$ in (\ref{eq3}) and  a real number $y_1\neq 0$. By Proposition \ref{spectralOf-Dm}(ii), there exist a eigenvector of $\textbf{L}_{S_{i_0}\rightarrow S_{i_0}}$ associated with the common eigenangle $\theta$, say $\textbf{y}^{(i_0)}$, such that $\textbf{y}^{(i_0)}|_{v}\neq 0(\forall v\in S_{i_0})$. For the same reason, there exist a eigenvector of $\textbf{L}_{S_{i_j}\rightarrow S_{i_j}}$ associated with $\theta$, say $\textbf{y}^{(i_j)}$, such that $\textbf{y}^{(i_j)}|_v\neq 0(\forall v\in S_{i_j})$ and

\begin{equation}\label{yi=-yi-1}
\textbf{y}^{(i_j)}|_{v_{i_j+1}}=-\textbf{y}^{(i_{j-1})}|_{v_{i_j-1}}(j=1,2,\cdots,k). \end{equation}

Set vector $\textbf{y}$ as $\textbf{y}|_{S_{i_j}}=\textbf{y}^{(i_j)}(j=0,1,\cdots,k)$ and $\textbf{y}|_{\overline{S}}=\textbf{0}$ . Armed with what we have proved above, we know that $\textbf{y}|_{v}\neq 0(\forall v\in S)$ and $\textbf{L}_{\{v_i\}\rightarrow S}\textbf{y}=0$(see (\ref{Lvi-S}) and ( \ref{yi=-yi-1})). Therefor,

\[\textbf{Ly}=\left[
\begin{array}{cc}
\textbf{L}_{S\rightarrow S}  &  \textbf{L}_{S\rightarrow \overline{S}}\\
\textbf{L}_{\overline{S}\rightarrow S}  &  \textbf{L}_{\overline{S}\rightarrow\overline{S}}
\end{array}
\right]
\left[
\begin{array}{c}
\textbf{y}|_{S}\\
\textbf{0}
\end{array}
\right]=\textbf{L}_{S\rightarrow S}\textbf{y}|_S.
\]
 This  means the vector $\textbf{y}$ is the eigenvector of $\textbf{L}$.

 Case 2 $M=0$.

 The proof is similar as Case 1 and trivial by noticing that $\textbf{L}_{S_{i_0}\rightarrow S_{i_0}}=\textbf{D}_m=\textbf{L}_{S_{i_k}\rightarrow S_{i_k}}$.
  \qed

If $S$ is a perfect critical vertex set of path $P_n$, then by Theorem \ref{theVertexNumberOfPath}, we have $n=m+(k-1)2m+m+k=k(2m+1)$. That is $n$ must be an odd. Therefor,
 a straightforward consequence of Theorem \ref{theVertexNumberOfPath} is that there exist a perfect critical vertex set of a path graph $P_n$ if and only if $n$ is odd. The following corollary follows straight from Theorem \ref{theVertexNumberOfPath}, which have been proved in \cite{parlangeli} by using different mathematical tools.

\begin{corollary}\label{n=2powerof2}
Let $n=2^{l_0}$ for some $l_0\in \mathbb{N}$, then the path $P_n$ is controllable with any vertex selected as leader, e.g. $P_n$ is omnicontrollable.   \qed
\end{corollary}

\subsection{Algorithm and Examples}

In fact, Theorem \ref{theVertexNumberOfPath} described all perfect critical vertex set of path graph $P_n$. Next, we provide  a method to locate the leader vertices. That is the following Algorithm I.

\begin{tabular}{l}
\hline
\textbf{Algorithm I}\\
\hline
1: input:  $n=2^{l_0}p_1^{l_1}p_2^{l_2}\cdots p_t^{l_t}$.\\

2: initialize: $N=\{p_1,p_2,\cdots,p_t\}$,  $F=\emptyset$ ,  $j=0$.\\

3: while  $N\neq \emptyset$, for some  $p\in N$   do\\

4:   \quad \quad $j=j+1$, $k=\frac{n}{p}, m=\frac{p-1}{2}$, $F_j=\emptyset$\\

5:   \quad \quad  for  $l=0: k-1$ do\\

6:   \quad \quad \quad \quad \quad  $i=(m+1)+l(2m+1)$\\

7:   \quad \quad \quad \quad \quad $F_j=F_j\bigcup\{v_i\}$\\

8: \quad \quad  end for\\

9:   \quad \quad  $N=N\setminus \{p\}$\\

10: \quad \,  $F=F\bigcup F_j $\\

11:   end while\\

12:  output:   $\overline{F}$\\
\hline
\end{tabular}

 For a path $P_n$, $\overline{F}$ obtained by Algorithm I is the set of leaders. That is $P_n$ is controllable with the vertex located in $\overline{F}$ and only with those vertices.

 For example, let $n=6$ is even. Since $n=2\times 3$ and only 3 is a odd prime factor, by Algorithm I, $p=3, k=\frac{n}{p}=2, m=\frac{p-1}{2}=1$, $F_1=\{v_i|i=(m+1)+l(2m+1),0\leq l \leq k-1\}=\{v_2,v_5\}$. So, any vertex in $\{v_1,v_3,v_4,v_6\}$ can be select as leader.

 Let $n=18$ is also even but  with 3 being a multiple factor.  By Algorithm I,
 there is only one follower vertex set needed to be considered. let $p=3, k=\frac{n}{p}=6, m=\frac{p-1}{2}=1$, $F_1=\{v_i|i=(m+1)+l(2m+1),0\leq l \leq k-1\}=\{v_2,v_5,v_8,v_{11},v_{14},v_{17}\}$. So, any vertex in $\overline{F_1}$ can be select as leader.

 In the case of multiple factor, Algorithm I is a  much more efficient algorithm to locate the leader vertex than the method provided in \cite{parlangeli}. What's more, Algorithm I can be easily applied to much lager path graph. For example, let $n=105$, since $n=3\times5\times7$, there are only three follower vertex set being calculate, that is

 $F_1=\{v_i|i=(1+1)+l\times(2\times1+1),0\leq l\leq 34 \}$,

 $F_2=\{v_i|i=(2+1)+l\times(2\times2+1),0\leq l\leq 20 \}$,

 $F_3=\{v_i|i=(3+1)+l\times(2\times3+1),0\leq l\leq 14 \}$.

 \hspace{-0.8cm}
$F=F_1\bigcup F_2 \bigcup F_3$ and leaders are located in the vertex set $\overline{F}$:

$\overline{F}=\{v_1,v_4 ,    v_6,     v_7,    v_9,    v_{10},    v_{12},    v_{15},
    v_{16},    v_{19},    v_{21},    v_{22},    v_{24},    v_{25},    v_{27},    v_{30},$

\hspace{1.1cm}$ v_{31},    v_{34},    v_{36},    v_{37},    v_{39},    v_{40},    v_{42},    v_{45},
     v_{46},    v_{49},    v_{51},    v_{52},    v_{54},    v_{55},    v_{57},$

\hspace{1.1cm} $v_{60},
      v_{61},    v_{64},    v_{66},    v_{67},    v_{69},    v_{70},    v_{72},    v_{75},
     v_{76},    v_{79},    v_{81},    v_{82},    v_{84},    v_{85},$

\hspace{1.1cm}  $ v_{87},    v_{90},
     v_{91},   v_{94},    v_{96},    v_{97},    v_{99},   v_{100},   v_{102},   v_{105}\}$.

 \hspace{-0.8cm}That is any one of and only of all these 56 vertices can be selected as leader.

\section{Minimal Perfect Critical Vertex Set of Graphs Based on Path}\label{section4}

Path graphs are simplest and basic graph structures. Some graphs can be constructed by adding paths. The minimal perfect critical vertex set of these graphs  will be investigated as what follows.

\begin{figure}[h]
\centering
\scriptsize
\setlength{\unitlength}{1mm}
\begin{picture}(200,80)(0,0)
\put(10,50){\dashbox{2}(40,25)}
\put(15,65){\circle{1.5}}
\put(30,55){\circle*{1.5}}
\put(24,62){\circle{1.5}}
\put(35,67){\circle{1.5}}
\put(45,62){\circle{1.5}}
\put(0,10){\circle{1.5}}
\put(22,70){$G$}
\put(12,28){\circle{1.5}}
\put(20,40){\circle{1.5}}
\put(0.5,10.5){\line(2,3){4}}
\put(5,17.5){\line(2,3){0.5}}
\put(6,19){\line(2,3){0.5}}
\put(7,20.5){\line(2,3){0.5}}
\put(8,22){\line(2,3){3.5}}
\put(12.5,28.8){\line(2,3){7}}
\put(20.4,40.5){\line(2,3){9.5}}

\put(15,1){$(a)\,\,\,G(v_0)+P_n$}
\put(-3,12){$v_1$}
\put(4,30){$v_{n-1}$}
\put(15,42){$v_n$}
\put(30,57){$v_0$}

\put(70,50){\dashbox{2}(40,25)}
\put(75,65){\circle{1.5}}
\put(90,55){\circle*{1.5}}
\put(84,62){\circle{1.5}}
\put(95,67){\circle{1.5}}
\put(105,62){\circle{1.5}}
\put(60,10){\circle{1.5}}
\put(82,70){$G$}
\put(72,28){\circle{1.5}}
\put(80,40){\circle{1.5}}
\put(60.5,10.5){\line(2,3){4}}
\put(65,17.5){\line(2,3){0.5}}
\put(66,19){\line(2,3){0.5}}
\put(67,20.5){\line(2,3){0.5}}
\put(68,22){\line(2,3){3.5}}
\put(72.5,28.8){\line(2,3){7}}
\put(80.3,40.5){\line(2,3){9.5}}

\put(70,1){$(b)\,\,\,G(v_0)+\{P_{n_1},P_{n_2},\cdots,P_{n_t}\}$}
\put(56,12){$v_1^{(1)}$}
\put(62,30){$v_{n_1-1}^{(1)}$}
\put(74,42){$v_{n_1}^{(1)}$}
\put(89,57){$v_0$}

\put(85,10){\circle{1.5}}
\put(85.1,10.5){\line(1,9){0.5}}
\put(85.7,16){\line(1,9){0.1}}
\put(86,17.8){\line(1,9){0.1}}
\put(86.2,19.6){\line(1,9){0.1}}
\put(86.4,21.2){\line(1,9){0.5}}
\put(87,26.2){\circle{1.5}}
\put(87.2,26.6){\line(1,9){1.3}}
\put(88.6,39){\circle{1.5}}
\put(88.4,39.6){\line(1,9){1.7}}

\put(78,12){$v_1^{(2)}$}
\put(89,26.2){$v_{n_2-1}^{(2)}$}
\put(90.6,39){$v_{n_2}^{(2)}$}

\put(100,20){\line(1,0){1}}
\put(102,20){\line(1,0){1}}
\put(104,20){\line(1,0){1}}

\put(125,10){\circle{1.5}}
\put(124.5,10.5){\line(-7,9){3}}
\put(121,14.9){\line(-7,9){0.5}}
\put(120,16.1){\line(-7,9){0.5}}
\put(119,17.3){\line(-7,9){0.5}}
\put(118,18.5){\line(-7,9){5.4}}
\put(112.3,26){\circle{1.5}}
\put(112,26.5){\line(-7,9){9.3}}
\put(102.3,39){\circle{1.5}}
\put(101.8,39.8){\line(-7,9){11.3}}

\put(125,12){$v_1^{(t)}$}
\put(114,26){$v_{n_t-1}^{(t)}$}
\put(104,39){$v_{n_t}^{(t)}$}

\end{picture}
\caption{Trees constructed by adding paths}
\label{treesByAddingPath}
\end{figure}
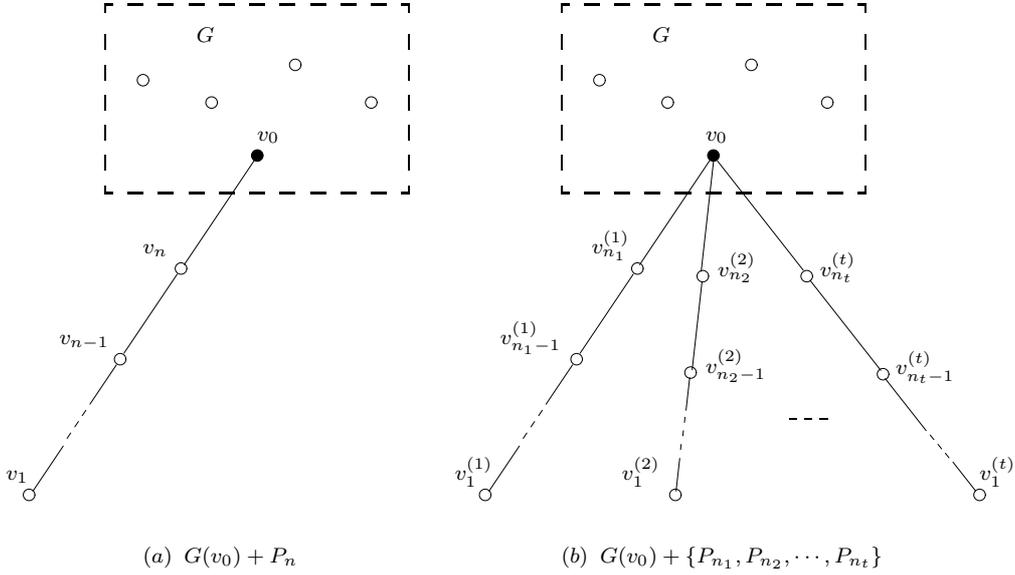

Let $G$ be a graph and $v_0\in V(G)$. We use $G(v_0)+\{P_{n_1},P_{n_2},\cdots,p_{n_t}\}$ to denote the graph by adding $P_{n_1},P_{n_2},\cdots,P_{n_t}$ to $G$ incident to $v_0$, as shown in fig.\ref{treesByAddingPath}.

%%%%%%%%%%%%%%  ͼ  5  %%%%%%%%%%%
%%%%%%%%%%%%%%  ͼ5 ½áÊø %%%%%%%%%%%%

Consider Lemma \ref{lemma1}, see fig.\ref{treesByAddingPath}(a), for any $S\subset \{v_1,v_2,\cdots,v_n\}$, we know that $S$ is not a perfect critical vertex set of $G(v_0)+P_n$. Therefor, we study the graphs $G(v_0)+\{P_{n_1},P_{n_2},\cdots,P_{n_t}\}$ and $t\geq 2$, see fig.\ref{treesByAddingPath}(b).

 Let $S$ be a minimal perfect critical vertex set of $G(v_0)+\{P_{n_1},P_{n_2},\cdots,P_{n_t}\}$. By Lemma \ref{lemma1}, the matrix $\textbf{L}_{S\rightarrow S}$ also has the form illustrated in (\ref{LStoS=blocked}) and the Lemma \ref{common-eigenvalue-of-Lsiktosik} is also hold. Also consider Lemma \ref{lemma1}, we know the last vertex $v_{n_l}^{(l)}$  belongs to $S$ when $S\bigcap V(P_{n_l})\neq\emptyset$. So,  from the what we have proved in section \ref{section3}, we know immediately that there exist exactly two pahts, say $P_{n_i},P_{n_j}$ , such that $S\bigcap V(P_{n_i})\neq\emptyset, S\bigcap V(P_{n_j})\neq\emptyset$ and $S\bigcap V(P_{n_k})=\emptyset(k\neq i,j)$. Further, we have the following theorem \ref{theorem5}.

\begin{theorem}\label{theorem5}
Let $G(v_0)+\{P_{n_1},P_{n_2},\cdots,P_{n_t}\}$ be the graph in fig.\ref{treesByAddingPath}(b) and $t\geq 2$.
Then there exist a minimal perfect critical vertex set $S$, where $S\subset V(P_{n_i})\bigcup V(P_{n_j})$, if and only if $2n_i+1$ and $2n_j+1$ have common divisor greater than 1.
\end{theorem}
\textbf{Proof}\quad Necessity: Noticing that if $S$ is a MPCVS of $G(v_0)+\{P_{n_1},P_{n_2},\cdots,P_{n_t}\}$, $S$ is a MPCVS of the path $P_{n_i}+v_0+P_{n_j}$. So, by Lemma \ref{|Dm|=|BM|} and Lemma \ref{lemma5}, we have $n_i=m+k_i(2m+1), n_j=m+k_j(2m+1)$ and $m>0$(by (\ref{m>1})). Hence $2n_i+1$ and $2n_j+1$ have a common divisor $2m+1$ greater than 1.

Sufficiency: Let $p>1$ is a common divisor of $2n_i+1$ and $2n_j+1$. $p$ is an odd.  Let $p=p_1^{l_1}p_2^{l_2}\cdots p_q^{l_q}$, where $p_i$ are primes. Set $m=\frac{p_1-1}{2}\geq 1$ and

\[S_1=V(P_{n_i})\backslash\{v_k^{(n_i)}|1\leq k \leq n_i, k=m+1(\texttt{mod}(2m+1))\}.\]

\[S_2=V(P_{n_j})\backslash\{v_k^{(n_j)}|1\leq k \leq n_j, k=m+1(\texttt{mod}(2m+1))\}.\]

Taking $S=S_1\bigcup S_2$, recall the sufficiency proof of Theorem \ref{theVertexNumberOfPath}, we know $S$ is a MPCVS.  \qed

Next, we provide some examples to illustrate how to use Theorem \ref{theorem5} to discover MPCVS of graphs constructed by adding paths.

%%%%%%%%%%%%%%  ͼ  6  %%%%%%%%%%%
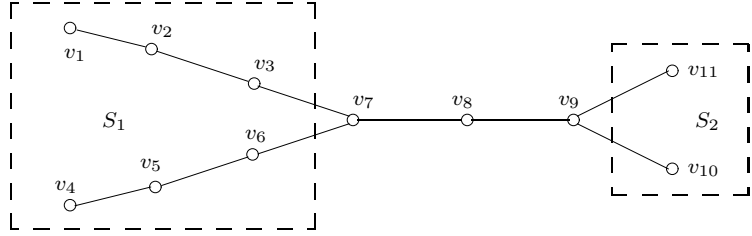
\begin{figure}[h]
\centering
\scriptsize
\setlength{\unitlength}{1mm}
\begin{picture}(100,35)
\put(0,4.5){\dashbox{2}(40,30)}
\put(45,19){\circle{1.5}}
\put(44.5,19.5){\line(-3,1){12}}
\put(32,23.8){\circle{1.5}}
\put(31.2,24.1){\line(-3,1){12}}
\put(18.5,28.4){\circle{1.5}}
\put(17.7,28.7){\line(-4,1){9}}
\put(7.8,31.2){\circle{1.5}}

\put(45,21){$v_7$}
\put(32,25.8){$v_3$}
\put(18.5,30.4){$v_2$}
\put(7,27.2){$v_1$}

\put(44.5,18.5){\line(-3,-1){12}}
\put(31.8,14.4){\circle{1.5}}
\put(31.2,14.1){\line(-3,-1){11.4}}
\put(19,10,1){\circle{1.5}}
\put(18.4,10.1){\line(-4,-1){10}}
\put(7.8,7.7){\circle{1.5}}

\put(30.8,16.4){$v_6$}
\put(17,12,6){$v_5$}
\put(5.8,9.7){$v_4$}

\put(60,19){\circle{1.5}}
\put(45.6,19){\line(1,0){13.5}}
\put(74,19){\circle{1.5}}
\put(60.6,19){\line(1,0){12.8}}

\put(58,21){$v_8$}
\put(72,21){$v_9$}

\put(79,9){\dashbox{2}(18,20)}

\put(74.5,19.5){\line(2,1){12}}
\put(87,25.5){\circle{1.5}}
\put(74.5,18.5){\line(2,-1){12}}
\put(87,12.5){\circle{1.5}}

\put(89,25){$v_{11}$}
\put(89,12){$v_{10}$}

\put(12,18){$S_1$}
\put(90,18){$S_2$}
\end{picture}
\caption{perfect critical vertex set of trees}
\label{ExampleForLeaderSetOfTrees}
\end{figure}
%%%%%%%%%%%%%%  ͼ6 ½áÊø %%%%%%%%%%%%

Example 1, fig.\ref{ExampleForLeaderSetOfTrees} comes from \cite{ZhijianAndHai}. By Theorem \ref{theorem5}, $S_1=\{v_1,v_2,\cdots,v_6\}$ and $S_2=\{v_{10},v_{11}\}$($S_2$ can also be discovered by Theorem \ref{theorem3})  are two MPCVSs. So, in order to make sure that the system is controllable, the minimum number of leader vertices is 2  and one of the leaders comes from $S_1$ and the other comes from $S_2$.

%%%%%%%%%%%%%%  ͼ  7  %%%%%%%%%%%
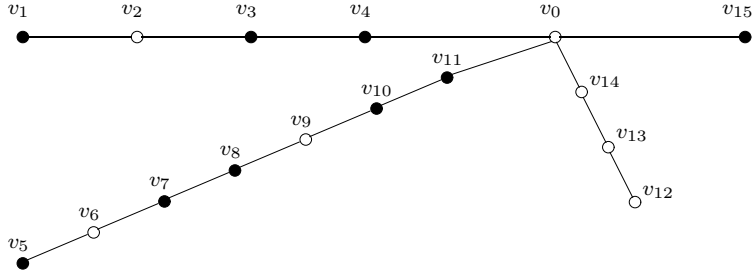
\begin{figure}[h]
\centering
\scriptsize
\setlength{\unitlength}{1mm}
\begin{picture}(120,40)
\put(10,35){\circle*{1.5}}
\put(25,35){\circle{1.5}}
\put(40,35){\circle*{1.5}}
\put(55,35){\circle*{1.5}}
\put(80,35){\circle{1.5}}
\put(105,35){\circle*{1.5}}

\put(8,38){$v_1$}
\put(23,38){$v_2$}
\put(38,38){$v_3$}
\put(53,38){$v_4$}
\put(78,38){$v_0$}
\put(102,38){$v_{15}$}

\put(10.6,35){\line(1,0){13.7}}
\put(25.6,35){\line(1,0){13.7}}
\put(40.6,35){\line(1,0){13.7}}
\put(55.6,35){\line(1,0){23.8}}
\put(80.6,35){\line(1,0){23.8}}

\put(10,5){\circle*{1.5}}
\put(10.5,5.4){\line(7,3){8}}
\put(19.3,9.1){\circle{1.5}}
\put(19.8,9.5){\line(7,3){8}}
\put(28.6,13.2){\circle*{1.5}}
\put(29.1,13.6){\line(7,3){8}}
\put(37.9,17.3){\circle*{1.5}}
\put(38.4,17.7){\line(7,3){8}}
\put(47.2,21.4){\circle{1.5}}
\put(47.7,21.8){\line(7,3){8}}
\put(56.5,25.5){\circle*{1.5}}
\put(57.0,25.9){\line(7,3){8}}
\put(65.8,29.6){\circle*{1.5}}
\put(66.3,30.0){\line(9,3){13.3}}

\put(8,7){$v_{5}$}
\put(17.3,11.1){$v_{6}$}
\put(26.6,15.2){$v_{7}$}
\put(35.9,19.3){$v_{8}$}
\put(45.2,23.4){$v_{9}$}
\put(54.5,27.5){$v_{10}$}
\put(63.8,31.6){$v_{11}$}

\put(80.2,34.3){\line(1,-2){3}}
\put(83.5,27.7){\circle{1.5}}
\put(83.7,27){\line(1,-2){3}}
\put(87,20.4){\circle{1.5}}
\put(87.2,19.7){\line(1,-2){3}}
\put(90.5,13.1){\circle{1.5}}

\put(84.5,28.7){$v_{14}$}
\put(88,21.6){$v_{13}$}
\put(91.5,14.1){$v_{12}$}

\end{picture}
 \caption{A Generalized Star}
 \label{figureForGeneralizedStar}
\end{figure}

Example 2. Generalized star is a kind of useful graphs constructed by paths. All MPCVSs of generalized star graph can be found by Theorem \ref{theorem5}.
 A \textit{star} with $n$ vertices is a graph consisting of one vertex $ v_0$ in the center and $n-1$ vertices adjacent to $v_0$ but not adjacent to each other. A \textit{generalized star} is the graph obtained from a star by replacing each edge by a path of arbitrary length. These paths are called legs, and generalized stars are also called \textit{spiders}. See fig.\ref{figureForGeneralizedStar}. Let $P_1=v_1v_2v_3v_4$, $P_2=v_5v_6\cdots v_{11}$, $P_3=v_{12}v_{13}v_{14}$, $P_4=v_{15}$. It is easily seen that $n_1=4, n_2=7, n_3=3, n_4=1$. By Theorem \ref{theorem5}, there exist MPCVS $S_1\subset V(P_1)\bigcup V(P_2)$, $S_2\subset V(P_1)\bigcup V(P_4)$, $S_3\subset V(P_2)\bigcup V(P_4)$. Further, by Algorithm I, we know that $S_1=\{v_1,v_3,v_4,v_5,v_7,v_8,v_{10},v_{11}\}$, $S_2=\{v_1,v_3,v_4,v_{15}\}$, $S_3=\{v_5,v_7,v_8,v_{10},v_{11},v_{15}\}$.
Arbitrarily selecting  two black vertices in fig.\ref{figureForGeneralizedStar}, say $v_i, v_j$, and $v_i, v_j$ belong to different paths, we know that the generalized star will be controllable with the leaders $\{v_i, v_j\}$. The minimum number of leaders is 2.

\section{Conclusion}\label{section5}
Neighbor-based controllability of undirected graph has received special attention in recent years. However, the understanding of roles of leaders, the minimum number of leaders, especially the leader location issues in the undirected graph is still largely unknown. A major effort in this paper is to provide a method to determine the leaders  directly from topology structures of undirected graph. These efforts also enlarge the understanding of the leader's role in undirected graph controllability.
To do this, we introduced the concept of critical vertex set,  perfect critical vertex set and minimal perfect critical vertex set. These concepts indicates that some vertices with special graphical characterization should be selected as leaders. Necessary and sufficient conditions are proposed to uncover some special minimal perfect critical vertex set. Theorem \ref{theorem3} described the graphical characterizations of minimal perfect critical 2 vertex set. Minimal perfect critical 3 vertex set do NOT exist was proved in Theorem \ref{theorem4}. Theorem \ref{theVertexNumberOfPath} completely described the MPCVS of path and Theorem \ref{theorem5} can be used to discover some special MPCVSs of graphs constructed by adding paths.

All these results clearly indicate where  leaders located, reveal the effect of topology structure on the controllability and promote a further study of controllability of undirected graph.

\end{document}